\documentclass[12pt]{amsart}
\usepackage[latin1]{inputenc}
\newtheorem{thm}{Theorem}
\newtheorem{cor}[thm]{Corollary}
\newtheorem{lem}[thm]{Lemma}

\theoremstyle{definition}
\newtheorem{defn}[thm]{Definition}

 \numberwithin{equation}{section}

\newcommand{\To}{\longrightarrow}
\begin{document}

\baselineskip=17pt

\title[]{Weakly countably determined spaces of high complexity}%
\author{Antonio Avil\'{e}s}

\thanks{The first author was supported by European Union Research Training Network PHD Marie
Curie, Intra-European Felloship MCEIF-CT2006-038768 and research
projects MTM2005-08379 (MEC and FEDER) and S\'{e}neca 00690/PI/04
(CARM)}

\address{University of Paris 7, Equipe de Logique
Math\'ematique, UFR de Mathématiques, 2 Place Jussieu, 75251
Paris, France} \email{aviles@logique.jussieu.fr, avileslo@um.es}

\begin{abstract}
We prove that there exist weakly countably determined spaces of
complexity higher than coanalytic. On the other hand, we also show
that coanalytic sets can be characterized by the existence of a
cofinal adequate family of closed sets. Therefore the Banach
spaces constructed by means of these families have at most
coanalytic complexity.
\end{abstract}

\maketitle

\section{Introduction}

We deal with the descriptive complexity of a Banach space $X$ with
respect to weak$^\ast$ compact subsets of the double dual
$X^{\ast\ast}$. The simplest Banach spaces in this sense are
reflexive spaces, which have a weakly compact ball and hence are
$\mathcal{K}_\sigma$ (that is, $\sigma$-compact) subsets of the
double dual. In the next level of complexity we find the class of
Banach spaces which are $\mathcal{K}_{\sigma\delta}$ subsets (that
is, countable intersection of $\mathcal{K}_\sigma$ sets) of their
double dual, which includes all weakly compactly generated (WCG)
spaces. Va\v{s}\'ak~\cite{Vasak} and Talagrand~\cite{WKA}
introduced, respectively, the following two further descriptive
classes:

\begin{defn}
A Banach space $X$ is called weakly countably determined (WCD) if
there exists a family $\{K_{s} : s\in\omega^{<\omega}\}$ of
weak$^\ast$ compact subsets of $X^{\ast\ast}$, and a set $A\subset
\omega^\omega$ such that $$X=\bigcup_{a\in
A}\bigcap_{n<\omega}K_{a|n}$$
\end{defn}

\begin{defn}
A Banach space $X$ is called weakly $\mathcal{K}$-analytic (WKA)
if there exists a family $\{K_{s} : s\in\omega^{<\omega}\}$ of
weak$^\ast$ compact subsets of $X^{\ast\ast}$, and an analytic set
$A\subset \omega^\omega$ such that $$X=\bigcup_{a\in
A}\bigcap_{n<\omega}K_{a|n}$$
\end{defn}

\renewcommand{\thefootnote}{\fnsymbol{footnote}}

Indeed, in the case of WKA spaces, the analytic set $A$ can be
chosen to be the whole Baire space $A=\omega^\omega$
\footnote[2]{The reader can try to check this directly as an
exercise, or else consider the characterizations of these concepts
in terms of \emph{uscos} \cite[pp. 117-118, 142]{FabianWA};
remember that any analytic set is a continuous image of
$\omega^\omega$} . Thus, the picture of the descriptive classes
already considered in these works is the following:
$$WCG\subset \mathcal{K}_{\sigma\delta} \subset WKA \subset WCD$$
One important problem in studying this hierarchy is the
``separation'' problem, that is, constructing examples showing
that the above inclusions are sharp. Talagrand used a technique to
solve this problem, the so called adequate families of
sets\footnote[3]{Talagrand attributes the concept of adequate
family to Roman Pol, cf. \cite[p. 417]{WKA}.}, which allowed him
to produce two examples:\\

\begin{itemize}
\item A Banach space which is $\mathcal{K}_{\sigma\delta}$ but not
weakly compactly generated~\cite{WKA}.

\item A weakly countably determined space which is not weakly
$\mathcal{K}$-analytic \cite{WCD}.\\
\end{itemize}

The remaining separation problem was left open by Talagrand and
has been recently solved by Argyros, Arvanitakis and
Mercourakis~\cite{ArgArvMer}, providing an example of a WKA space
which is not a $\mathcal{K}_{\sigma\delta}$ space. They construct
their example by a technique different from Talagrand's adequate
families, by using the so called Reznichenko families of trees.
Indeed, they prove a result which shows that it is impossible to
produce a $\mathcal{K}_{\sigma\delta}$ non WKA space using
adequate families, which explains Talagrand's failure in solving
this question. Argyros, Arvanitakis and Mercourakis have also
succeeded in showing that it is not only that the classes
$\mathcal{K}_{\sigma\delta}$ and $WKA$ can be separated, but
indeed there is a whole Borel hierarchy of spaces between them,
$$\mathcal{K}_{\sigma\delta}\subset
\mathcal{K}_{\sigma\delta\sigma\delta}\subset\cdots\subset WKA.$$
In this note, we shall focus on higher levels of this hierarchy,
on the gap between WKA and WCD. If one looks at Talagrand's
example~\cite{WCD} separating these two classes, one realizes that
the set $A\subset\omega^\omega$ which witnesses that it is WCD is
indeed a complete coanalytic set. We propose the following
definition:

\begin{defn}
Let $\mathcal{C}$ be a class of separable metrizable spaces. A
Banach space is said to be \emph{weakly $\mathcal{C}$-determined}
if there exists a family $\{K_{s} : s\in\omega^{<\omega}\}$ of
weak$^\ast$ compact subsets of $X^{\ast\ast}$, and a set
$A\in\mathcal{C}$, $A\subset\omega^\omega$ such that $$X=\bigcup_{a\in
A}\bigcap_{n<\omega}K_{a|n}.$$
\end{defn}

In this language, Talagrand's example~\cite{WCD} is a weakly
$\Pi_1^1$-determined space which is not weakly
$\Sigma_1^1$-determined (The symbols $\Pi_1^1$ and $\Sigma_1^1$
represent the classes of coanalytic and analytic
sets\footnote[2]{It is usual to consider the notion that a subset
$A$ of a Polish space $X$ is Borel, analytic, coanalytic,
$\Sigma^1_n$, $\Pi^1_n$, etc. However all these properties are
intrinsic topological properties of $A$ which do not depend on the
Polish superspace \cite{Kechrisbook}. Thus we talk about separable
metrizable spaces which are Borel, analytic and so on.} in the
logical notation, cf.\cite{Kechrisbook}). The natural question
arises: Are there WCD spaces of higher complexity? Namely, are
there WCD spaces that are not weakly $\Pi_1^1$-determined? The two
main results of this note are motivated by this problem and draw a
similar picture as the one obtained by Argyros, Arvanitakis and
Mercourakis in the lower level of the hierarchy.\\

In the first part of our work we analyze the technique of adequate
families that Talagrand used for his two examples. We simplify
this construction and we show that the right framework for it is
that of coanalytic sets. Again, although for different reasons
than in the $\mathcal{K}_{\sigma\delta}$ problem, our
Theorem~\ref{adequatecoanalytic} shows that Talagrand's technique
cannot produce WCD spaces of higher complexity than coanalytic.
This result is indeed an intrinsic topological characterization of
coanalytic sets which may have its independent interest.\\

Our second result states that the technique of Reznichenko
families of trees developed in~\cite{ArgArvMer} does allow to give
a positive answer to our question: There are WCD spaces of
complexity higher than coanalytic, indeed there are WCD spaces of
\emph{arbitrarily high complexity}, in a sense that will be made
precise. In particular, all projective classes can be separated:

\begin{thm}\label{projectiveclasses}
For every $n\geq 1$ there exists a Banach space which is weakly
$\Sigma^1_{n+1}$-determined but not weakly
$\Sigma^1_n$-determined.
\end{thm}

\section{General facts about WCD spaces}

\begin{defn} A class $\mathcal{C}$ of separable metrizable spaces
will be called \emph{nice} if it is closed under the following
operations:

\begin{itemize}
\item closed subspaces, \item continuous images, \item countable
products, \item Wadge reduction, that is, if $f:A\To B$ is a
continuous function between Polish spaces, $C\subset B$,
$C\in\mathcal{C}$, then $f^{-1}(C)\in \mathcal{C}$.
\end{itemize}
\end{defn}

\begin{defn}
Let $\Gamma$ be an uncountable set, let
$K\subset\mathbb{R}^\Gamma$ be a compact subset, and $D$ a
separable metrizable space. A mapping $f:\Gamma\To D$ is called a
\emph{determining function} if for every $x\in K$, for every
compact subset $C\subset D$, and every $\varepsilon>0$ we have
that $\{\gamma\in f^{-1}(C) : |x_\gamma|>\varepsilon\}$ is finite.
\end{defn}

We leave to the reader to check that the fact that $f$ like above is a determining function is equivalent to any the two following statements:\\

\begin{itemize}
\item[(i)] For every $x\in K$ and for every $\varepsilon>0$ there
is a neighborhood $V$ of $x$ in $K$ such that $\{\gamma\in
f^{-1}(V) : |x_\gamma|>\varepsilon\}$ is finite. \item[(ii)] For
every $x\in K$ and $\varepsilon>0$, the restriction of $f$ to the
set $\{\gamma\in\Gamma : |x_\gamma|>\varepsilon\}$ is
finite-to-one and has a closed and discrete range.\\
\end{itemize}

All the Banach spaces that we consider in this note are spaces of
continuous functions $C(K)$. The following theorem provides a
useful criterion to identify when a space $C(K)$ is a weakly
$\mathcal{C}$-determined space. The history of this result goes
back to \cite{WKA}, \cite{Farmaki}, \cite{ArgFar}, \cite{c1Mer}
and \cite{WCGrelatives}. Originally it has been stated for WKA or
WCD spaces, but it holds for any nice class $\mathcal{C}$.

\begin{thm}\label{determiningfunction}
Let $\mathcal{C}$ be a nice class and let $K\subset
\mathbb{R}^\Gamma$ be a compact set such that every $x\in K$ has
countable support, that is, $|\{\gamma\in\Gamma : x_\gamma\neq
0\}|\leq\omega$. Then the following are equivalent:
\begin{enumerate}
\item $C(K)$ is a weakly $\mathcal{C}$-determined space. \item
There exists $D\in\mathcal{C}$ and a determining function
$f:\Gamma\To D$.
\end{enumerate}
\end{thm}

{\tt Proof}: First of all, one should notice that in part (2) of
the theorem, $D$ can be taken to be a subset of $\omega^\omega$.
The reason is that, because we considered $\mathcal{C}$ to be
closed under Wadge reductions, for every $D_1\in\mathcal{C}$ there
exists $D_2\subset \omega^\omega$, $D_2\in\mathcal{C}$ such that
$D_2$ maps continuously onto $D_1$, and we have the following
fact:

\begin{lem}\label{determiningcomposition} If $f:\Gamma \To D_2$ is a determining function and $g:D_1\To D_2$ is a continuous surjection, then there exists a determining function $f':\Gamma \To
D_1$.
\end{lem}

{\tt Proof of the lemma}: Take $s:D_2\To D_1$ any mapping (not
necessarily continuous) such that $gs=1_{D_2}$. We prove that
$sf:\Gamma \To D_1$ is a determining function. Suppose it was not.
Then there would exist $x\in K$, $\varepsilon>0$ and $C\subset
D_1$ compact such that $$\{\gamma\in (sf)^{-1}(C) :
|x_\gamma|>\varepsilon\}$$ is infinite. But
$$\{\gamma\in (sf)^{-1}(C) : |x_\gamma|>\varepsilon\}\subset \{\gamma\in f^{-1}(g(C)) : |x_\gamma|>\varepsilon\}.$$
This contradicts that $f$ is a determining function. The lemma is
proved.$\qed$\\

After this observation, the statement of the theorem is the same
as \cite[Theorem 10(c)]{WCGrelatives} after changing ``$K$ is
Talagrand compact'' by ``$C(K)\in \mathcal{C}$', changing
``$\Gamma =
\bigcup_{\sigma\in\omega^\omega}\bigcap_{j=1}^\infty\Gamma_{\sigma|j}$''
by ``$\Gamma = \bigcup_{\sigma\in
D}\bigcap_{j=1}^\infty\Gamma_{\sigma|j}$ for some
$D\in\mathcal{C}$'', and changing ``$\forall
\sigma\in\omega^\omega$'' by ``$\forall \sigma\in D$''. It is now
a long but straightforward exercise that if one goes along the
proofs of \cite[Theorem 4]{WCGrelatives} and \cite[Theorem
10(c)]{WCGrelatives} and makes the obvious modifications, no
obstacle appears and a proof of Theorem~\ref{determiningfunction}
is obtained.$\qed$\\

The compact spaces $K$ for which $C(K)$ is WCG, WKA and WCD are
called Eberlein, Talagrand and Gul'ko compact respectively. We
shall call $\mathcal{C}$-Gul'ko compact those compact spaces $K$
for which $C(K)$ is weakly $\mathcal{C}$-determined.

\section{Adequate families on coanalytic sets}

A family $\mathcal{A}$ of subsets of a set $X$ is called an
\emph{adequate} family of sets if it satisfies the two following
properties:

\begin{itemize}
\item If $A\in\mathcal{A}$ and $B\subset A$, then
$B\in\mathcal{A}$.

 \item If $B$ is a subset of $X$ such that all finite
 subsets of $B$ belong to $\mathcal{A}$, then $B\in\mathcal{A}$.

\end{itemize}

In other words, to state that $\mathcal{A}$ is an adequate family
is equivalent to state that a subset $B\subset X$ belongs to
$\mathcal{A}$ if and only if every finite subset of $B$ belongs to
$\mathcal{A}$. Every adequate family of subsets of $X$ can be
naturally viewed as a closed subset of the product $\{0,1\}^X$ and
hence, is a compact Hausdorff space.\\

The interesting case for us occurs when $X$ is a separable
metrizable space and $\mathcal{A}$ is an adequate family which
consists of some closed subsets of $X$ (indeed closed and
discrete, since the family is hereditary) because then we get a
weakly countably determined space:

\begin{thm}\label{adequateisCdetermined}
Let $\mathcal{C}$ be a nice class, let $X\in\mathcal{C}$ and
$\mathcal{A}$ be an adequate family of closed subsets of $X$. Then
$C(\mathcal{A})$ is a weakly $\mathcal{C}$-determined Banach
space.
\end{thm}

This follows immediately from Theorem~\ref{determiningfunction}
just taking the identity $f:X\To X$ as a determining function.
Talagrand's example from~\cite{WKA} is an adequate family of
closed subsets of $X=\omega^\omega$ and the one from~\cite{WCD} is
an adequate family of closed subsets of $X=WF$, the set of well
founded trees on $\omega^{<\omega}$, the standard complete
coanalytic set. The fact that the first Banach space is WKA and
the second Banach space is WCD follows immediately from the above
theorem. But the negative part, that they are not WCG and WKA
spaces respectively needs further arguments and relies on the fact
that these adequate families are taken to be \emph{big enough} (of
course, not any adequate family of closed sets would work). We
have isolated the property of these adequate families which makes
them be as complicated as their underlying set.

\begin{defn} We say that an adequate family $\mathcal{A}$ of closed
subsets of a topological space $X$ is \emph{cofinal} if for every
infinite closed and discrete subset $B$ of $X$ there exists an
infinite subset $A\subset B$ such that $A\in\mathcal{A}$.
\end{defn}

The following theorem is the main result of this section. The
implication $2\Rightarrow 1$ constitutes a generalization and at
the same simplification of Talagrand's construction
from~\cite{WCD} (in particular we are avoiding any manupulation
with trees, using instead the easier and more general coanalytic
structure). The converse $1\Rightarrow 2$ establishes the
impossibility of arranging this construction out of the framework
of coanalytic sets.

\begin{thm}\label{adequatecoanalytic}
For a separable metrizable space $X$ the following are equivalent:
\begin{enumerate}
\item There exists a cofinal adequate family of closed subsets of
$X$. \item $X$ is coanalytic.
\end{enumerate}
\end{thm}

{\tt Proof}: Let $(K,d)$ be a compact metric space which contains $X$ as
a dense set, $K=\overline{X}$, and let $Y=K\setminus X$. We denote
by $\mathcal{M}$ the space of all strictly increasing
sequences of positive integers, which is homeomorphic to the Baire space $\mathbb{N}^\mathbb{N}$.\\

[$1\Rightarrow 2$] Let $\mathcal{A}$ be a cofinal adequate family
of closed subsets of $X$ and let $\{a_n : n=1,2,\ldots\}$ be an
enumeration of a dense subset of $X$. We consider the set
$$C = \left\{(y,\sigma)\in K\times\mathcal{M} : \{a_{\sigma_1},a_{\sigma_2},\ldots\}\in\mathcal{A}\text{ and }
d(y,a_{\sigma_i})\leq \frac{1}{i}\text{ for all }i\geq
1\right\}.$$ In order to prove that $X$ is coanalytic, we check
that $Y$ is analytic by showing that $C$ is a closed subset of
$K\times\mathcal{M}$ and $Y = \{y\in K : \exists
\sigma\in\mathcal{M} : (y,\sigma)\in C\}$. If we pick
$(y,\sigma)\in K\times\mathcal{M}\setminus C$ then either
$\{a_{\sigma_1},a_{\sigma_2},\ldots\}\not\in\mathcal{A}$ or there
exists $i\in\mathbb{N}$ such that $d(y,a_{\sigma_i})>
\frac{1}{i}$. In the first case, since $\mathcal{A}$ is an
adequate family, there exists $j\in\mathbb{N}$ such that
$\{a_{\sigma_1},a_{\sigma_2},\ldots,a_{\sigma_j}\}\not\in\mathcal{A}$
and then $\{(z,\tau)\in K\times\mathcal{M} : \tau_r = \sigma_r
\forall r\leq j\}$ is a neighborhood of $(y,\sigma)$ which does
not intersect $C$. In the second case, there exists a neighborhood
$U$ of $y$ such that $d(z,a_{\sigma_i})> \frac{1}{i}$ for all
$z\in U$ and then $\{(z,\tau) : z\in U\text{ and }
\tau_i=\sigma_i\}$ is a neighborhood of $(y,\sigma)$ which does
not intersect $C$. This proves that $C$ is a closed set.\\

 We show now that $Y = \{y\in K : \exists \sigma\in\mathcal{M} :
(y,\sigma)\in C\}$. Let us fix $y\in Y$. The sequence
$\{a_1,a_2,\ldots\}$ is a dense subset of $X$ which is moreover
dense in $K$, so there exists a subsequence
$\{a_{n_1},a_{n_2},\ldots\}$ which converges to $y$. Since
$y\not\in X$, the set $\{a_{n_1},a_{n_2},\ldots\}$ is a closed and
discrete subset of $X$, hence, since $\mathcal{A}$ is a
\emph{cofinal} adequate family in $X$, this sequence has a
subsequence $\{a_{m_1},a_{m_2},\ldots\}\in\mathcal{A}$ which still
converges to $y$. We can pass still to a further subsequence
$\{a_{k_1},a_{k_2},\ldots\}\in\mathcal{A}$ such that
$d(a_{k_i},y)\leq\frac{1}{i}$. If we call $\sigma =
(k_1,k_2,\ldots)$, we found that $(y,\sigma)\in C$. Conversely,
let us suppose now that we have a pair $(y,\sigma)\in C$. Then the
sequence $\{a_{\sigma_1},a_{\sigma_2},\ldots\}$, being a member of
the adequate family $\mathcal{A}$, constitutes a closed and
discrete subset of $X$, but at the same time the sequence
converges to
$y$, so we deduce that $y\not\in X$, and hence $y\in Y$.\\

 [$2\Rightarrow 1$] Suppose $C\subset K\times\mathcal{M}$
is a closed set such that $Y = \{x\in K : \exists
\sigma\in\mathcal{M} \text{ with } (x,\sigma)\in C\}$. Let $\prec$
be a well order on $X$ (the use of the axiom of choice here is not
essential, but it deletes a number of technicalities in the
proof). We define the cofinal adequate family $\mathcal{A}$ in the
following way. A finite set belongs to $\mathcal{A}$ if and only
if it is of the form $\{x_1\prec\cdots\prec x_n\}$ and there
exists $(y,\sigma)\in C$ such that
$d(y,x_i)\leq\frac{1}{\sigma_i}$ for every $i=1,\ldots,n$ (notice
that this is an hereditary condition, if a finite set satisfies
it, then so does every subset). An infinite set belongs to
$\mathcal{A}$ if and only if every finite
subset belongs to $\mathcal{A}$.\\

First, we show that every infinite set $A\in\mathcal{A}$ is a
closed and discrete subset of $X$. Otherwise, there would exist a
sequence $\{x_1\prec x_2\prec\cdots \}\subset A$ which converges
to a point $x\in X$ with $x\neq x_i$ for all $i$. Since
$A\in\mathcal{A}$, for every $n$ there exists $(y^n,\sigma^n)\in
C$ such that $d(y^n,x_i)\leq\frac{1}{\sigma^n_i}$ for every $i\leq
n$. Notice that the sequence $(y^n)$ also converges to $x$ because
$d(y^n,x_n)\leq \frac{1}{\sigma^n_n}\leq\frac{1}{n}$ (recall that
all sequences in $\mathcal{M}\subset\omega^\omega$ are
strictly increasing). Observe also that for every $i\in\omega$
the transversal sequence $\{\sigma^n_i : n=1,2,\ldots\}$ is
eventually constant to a value that we call $\sigma^\infty_i$
(otherwise the sequence $(y^n)$ would converge to $x_i$, since
$d(y^n,x_i)\leq \frac{1}{\sigma_i^n}$ for every $n\geq i$). The
sequence $\{\sigma^1,\sigma^2,\ldots\}$ converges in $\mathcal{M}$
to the sequence $\sigma^\infty =
(\sigma^\infty_1,\sigma^\infty_2,\ldots)$, and also
$(y^n,\sigma^n)$ converges to $(x,\sigma^\infty)$ and since $C$ is
a closed set, we find that $(x,\sigma^\infty)\in C$ which
contradicts that $x\in X$.\\

It remains to show that the adequate family $\mathcal{A}$ has the
property of being cofinal in $X$. We take $B$ an infinite closed
and discrete subset of $X$. Viewing $B$ as a subset of the compact
space $K$, we know that there exists a sequence $\{x_1\prec
x_2\prec\cdots\}\subset B$ which converges to a point $y\in K$.
Since $B$ is closed and discrete in $X$, it must be the case that
$y\in Y$. Therefore, we can pick $\sigma\in\mathcal{M}$ such that
$(y,\sigma)\in C$ and then we can find a subsequence
$\{x_{n_1}\prec x_{n_2}\prec\cdots\}$ such that $d(y,x_{n_i})\leq
\frac{1}{\sigma_i}$ for every $i$. This subsequence is an element
of $\mathcal{A}$, since every finite cut of the sequence satisfies
the definition of the family $\mathcal{A}$ with the same witness
$(y,\sigma)\in C$.$\qed$\\

We devote the rest of the section to check that for a cofinal
adequate family $\mathcal{A}$ of closed subsets of a separable
metrizable space $X$, the complexity of $C(\mathcal{A})$ is the
same as the complexity of $\mathcal{A}$ and not lower. We mention
that \v{C}i\v{z}ek and Fabian~\cite{CizFab} already realized that,
by transferring the original examples of Talagrand, given any
0-dimensional complete metrizable space $X$, then for every
coanalytic non-Borel subset $Y\subset X$ there is an adequate
family of subsets of $Y$ sot that the corresponding compact is
Gul'ko but not Talagrand, and that for every Borel
non-$\sigma$-compact subset $Y\subset X$ there is an adequate
family of subsets of $Y$ so that the corresponding compact space
is Talagrand but not Eberlein. They also gave a simpler approach
in checking the negative part in the first
kind of examples, which we shall follow.\\

 For a
family $\mathcal{A}$ of subets of a set $X$ and a subset $Z\subset
X$, we denote $\mathcal{A}|_Z = \{A\cap Z : A\in\mathcal{A}\}$,
the restriction of the family $\mathcal{A}$ to the set $Z$.
 When $\mathcal{A}$ is an adequate family we can write $\mathcal{A}|_Z = \{A\in\mathcal{A} : A\subset
 Z\}$.\\

\begin{thm}\label{adequateEberlein}
Let $X$ be a coanalytic space, $\mathcal{A}$ be a cofinal adequate
family of closed subsets of $X$ and let $Z$ be a subset of $X$.
Then $\mathcal{A}|_Z$ is an Eberlein compact if and only if $Z$ is
contained in some $\sigma$-compact subset of $X$.
\end{thm}

\begin{cor}
Let $X$ be a Borel non $\sigma$-compact space, $\mathcal{A}$ a
cofinal adequate family of closed subsets of $X$ and $Z\subset X$
any subset not contained in any $\sigma$-compact subset of $X$.
Then $\mathcal{A}|_Z$ is a Talagrand non Eberlein compact space.
\end{cor}

We notice that it follows from~\cite[Theorem 1.4]{ArgArvMer} that
if an adequate family $\mathcal{A}$ is Talagrand compact, then
indeed $C(\mathcal{A})$ is a $\mathcal{K}_{\sigma\delta}$ space.\\

{\tt Proof of Theorem~\ref{adequateEberlein}}: Assume $\mathcal{A}|_Z$ is an Eberlein compact. $Z$ being metrisable and separable, every set in $\mathcal{A}|_Z$ is at most countable. Then, there is a decompostion $Z=\bigcup_{n<\omega}Z_n$ such that for every $n<\omega$ the family $\mathcal{A}|_{Z_n}$ contains finite sets only \cite[Theorem 4.3.2]{FabianWA}. Fix $n<\omega$; we show that $Z_n$ is a relatively compact subset of $X$. Let $(z_m)_{m<\omega}$ be a one-to-one sequence in $Z_n$, and suppose for contradiction that it contains no subsequence convergent in $X$. Then, it must contain a subsequence which is closed and discrete in $X$. From the cofinality of $\mathcal{A}$ this subsequence contains an infinite subset $A\in\mathcal{A}$. Hence $A\cap Z_n$ is infinite, a contradiction.\\

Conversely, suppose that $Z\subset \bigcup_{n<\omega}K_n$ where each $K_n$ is a compact subset of $X$. Fix any $A\in\mathcal{A}$ and any $n<\omega$. We claim that the set $A\cap Z\cap K_n$ is finite. If not, because $A$ is closed and $K_n$ is compact, this set would contain a sequence convergent to some $x\in K_n\cap A$ which is in contradiction with the discreteness of $A$. Having the claim proved, we obtain that $\mathcal{A}|_Z$ is an Eberlein compact \cite[Theorem 4.3.2]{FabianWA}.$\qed$\\

\begin{thm}\label{adequateTalagrand}
Let $X$ be a coanalytic space, $\mathcal{A}$ be a cofinal adequate
family of closed subsets of $X$ and let $Z$ be a subset of $X$.
Then $\mathcal{A}|_Z$ is Talagrand compact if and only if $Z$ is
contained in some Borel subspace of $X$.
\end{thm}

\begin{cor}
Let $X$ be a coanalytic non Borel space, $\mathcal{A}$ a cofinal
adequate family of closed subsets of $X$ and $Z$ any subset of $X$
not contained in any Borel subspace of $X$. Then $\mathcal{A}|_Z$
is a Gul'ko compact (indeed $\Pi_1^1$-Gul'ko) which is not
Talagrand compact.
\end{cor}

{\tt Proof of Theorem~\ref{adequateTalagrand}}: If $Z$ is contained in
some Borel space $B\subset X$, then $\mathcal{A}|_Z$ can be viewed
as an adequate family of closed subsets of $B$, and then it
follows from Theorem~\ref{adequateisCdetermined} that
$C(\mathcal{A}|_Z)$ is WKA.\\

Now suppose that $\mathcal{A}|_Z$ is a Talagrand compact. Then by
Theorem~\ref{determiningfunction} there is a determining function $f:Z\To A$ to an analytic set $A$, indeed there is a
determining function $\psi:Z\To\omega^\omega$, that we can
get by composing $f$ with a selection for a continuous surjection
$\omega^\omega\To A$.  We consider as usual $K$ a compact
metric space with $\overline{X}=K$. For a finite sequence of
natural numbers $(k_1,\ldots,k_n)$ we denote $[k_1,\ldots,k_n] =
\{\tau\in\omega^\omega : \forall i\leq n\ \tau_i=k_i\}$. We
claim that
$$Z\subset
\bigcup_{\sigma\in\omega^\omega}\bigcap_{n\in\omega}\overline{\psi^{-1}[\sigma_1,\ldots,\sigma_n]}\subset
X,$$ where the closures are taken inside $K$. The first inclusion
is clear since
$z\in\bigcap_{n\in\omega}\overline{\psi^{-1}[(\psi(z)_1,\ldots,\psi(z)_n)]}$
for every $z\in Z$. For the second inclusion, suppose by
contradiction that for some $\sigma\in\omega^\omega$ we
have $y\in
\bigcap_{n\in\omega}\overline{\psi^{-1}[\sigma_1,\ldots,\sigma_n]}\setminus
X$. Then we can find a sequence of elements $y_n\in
\psi^{-1}[\sigma_1,\ldots,\sigma_n]$ which converges to $y$. Since
$y\not\in X$, the set $\{y_n : n\in\omega\}$ is an infinite
closed and discrete subset of $X$, so by cofinality we find a
subsequence $a=\{y_{n_k} : k\in\omega\}\in\mathcal{A}$. Then,
the image of the support of $a\in\mathcal{A}|_Z$ under $\psi$ is a
convergent sequence in $\omega^\omega$, which contradicts
that $\psi$ is a determining function.\\

We found that $Z$ is contained in a subset of $X$ which is a
Souslin operation of closed subsets of $K$, hence analytic. Since
$X$ is coanalytic, by the separation theorem (every two disjoint
analytic sets in a Polish space can be separated by disjoint
larger Borel sets) we deduce that $Z$ is contained in a Borel
subspace of $X$.$\qed$\\

\section{Gul'ko compact spaces of higher complexity}\label{SectionReznichenko}

We recall know the notion of Reznichenko family of trees
associated to a hereditary family of sets and the corresponding
compact space, which have been introduced and studied in
\cite{ArgArvMer}. In what follows, by a \emph{tree} we mean a set
$T$ endowed with a partial order relation $\leq$ such that (1) for
every $t\in T$ the set $\{s\in T : s< t\}$ is well ordered, and
(2) $T$ has a $\prec$-minimum, called the \emph{root} of $T$. An
element of the tree $t\in T$ is called a \emph{node} of $T$. An
immediate successor of $t\in T$ is a node $s< t$ for which there
is no further node $r$ with $t< r< s$. For an ordinal $\alpha$,
the $\alpha$-th level of the tree $T$ is the set of all $t\in T$
such that $\{s:s< t\}$ has order type $\alpha$. The height of a
tree is the first ordinal $\alpha$ for which the $\alpha$-th level
is empty. A subset $S$ of a tree $(T,\leq)$ is a \emph{segment} if
(1) every two elements of $S$ are comparable in the order $\leq$,
and (2) if $t,s\in S$, $r\in T$ and $t\leq r\leq  s$ then $r\in
S$. A segment $S$ is \emph{initial} if it contains the root of
$T$.\\

Let $A$ be a set of cardinality at most $\mathfrak c$ and
$\mathcal{F}$ a ``hereditary'' family of subsets of $A$
(hereditary means that if $B\in\mathcal{F}$ and $C\subset B$, then
$C\in\mathcal{F}$). An $(A,\mathcal{F})$-\emph{Reznichenko family
of trees} is a family of trees $\{T_a : a\in A\}$ indexed by the
set $A$ with the following properties:

\begin{enumerate}
\item For every $a\in A$, $T_a$ is a tree of height $\omega$, and
in which every node has $\mathfrak c$ many immediate successors
(in particular, $T_a$ has cardinality $\mathfrak c$). \item
$T_a\cap A =\{a\}$ and $a$ is the root of $T_a$. \item For every
$t\in\bigcup_{a\in A}T_a$, we have that $\{a\in A : t\in
T_a\}\in\mathcal{F}$.\item For every $a,b\in A$, $a\neq b$ and
every $S\subset T_a$ and $S'\subset T_b$ segments of the trees
$T_a$ and $T_b$ respectively, we have that $|S\cap S'|\leq 1$.
\item For every $B\in\mathcal{F}$ and for every disjoint family
$\{S_b : b\in B\}$ where $S_b$ is a finite initial segment of the
tree $T_b$ for every $b\in B$, there exist $\mathfrak c$ many
elements $t$ with the property that $t$ is an immediate successor
of the segment $S_b$ in the tree $T_b$ simultaneously for all
$b\in B$.
\end{enumerate}

In this context we put $\mathcal{T} = \bigcup_{a\in A}T_a$ and let
$R[\mathcal{F}]\subset 2^\mathcal{T}$ denote the family of all
segments of all the trees $T_a$, which can be easily checked to be
a compact family. It is shown in~\cite{ArgArvMer} that there does
exist an $(A,\mathcal{F})$-Reznichenko family of trees for any
given set $A$ of cardinality $\mathfrak c$ and any
hereditary family $\mathcal{F}$ of subsets of $A$.\\

Recall that $\omega^{<\omega}$ is the set of finite sequences of
natural numbers, ordered in the following way:
$(s_i)_{i<n}<(t_i)_{i<n}$ if $n\leq m$ and $s_i=t_i$ for $i<n$. In
order to avoid confusion with the concept of tree introduced
before, we define a \emph{tree on $\omega$} to be a subset
$T\subset \omega^{<\omega}$ such that if $a\in T$ and $b<a$ then
$b\in T$. We denote by $Tr\subset 2^{\omega^{<\omega}}$ the family
of all trees on $\omega$; this is a compact family and is viewed
as a compact metrizable space. A branch of $T\in Tr$ is an
infinite sequence $a\in\omega^\omega$ such that $a|n\in
T$ for all $n$.\\

We fix $A\subset\mathcal{N}$ to be a subset of the Baire space
$\mathcal{N}=\omega^\omega$. For every
$s\in\omega^{<\omega}$ we shall denote $W_s = \{a\in A : s\prec
a\}$ where $s\prec a$ means that if $s=(s_i)_{i<n}$, then $s_i =
a_i$ for all $i<n$. These sets constitute a basis for the topology
of $A$. Also, we shall denote by $wf(A)\subset Tr$ the family of
all trees on $\omega$ none of whose branches are elements of $A$.
The following theorem asserts that in this context, a compact
space $R[\mathcal{F}]$ constructed as above from a hereditary
family of closed and discrete subsets of $A$ is always a Gul'ko
compact, whose complexity is bounded by the
complexity of the set $wf(A)$. This is nothing else than a more informative restatement of some lemmas from~\cite{ArgArvMer}. Nevertheless, we found convenient to include a complete proof.\\

\begin{thm}\label{boundcomplexity}
Let $A\subset \mathcal{N}$, let $\mathcal{T} = \bigcup_{a\in
A}T_a$, and let $\mathcal{F}$ be an hereditary family of closed
and discrete subsets of $A$ and let $R[\mathcal{F}]\subset
2^\mathcal{T}$ be the compact set coming from an
$(A,\mathcal{F})$-Reznichenko family of trees. Then, there exists
a determining function $f:\mathcal{T}\To wf(A)\times
\omega^{(\omega^{<\omega})}$.
\end{thm}

{\tt Proof}: Let $t\in\mathcal{T}=\bigcup_{a\in A}T_a$, and let $B(t) =
\{a\in A : t\in T_a\}$ which is a set from $\mathcal{F}$ and hence
closed and discrete in $A$. We define

$$f_1(t) = \{s\in\omega^{<\omega} : |W_s\cap B(t)|>1\}.$$

Clearly, $f_1(t)\in wf(A)$ because if $a\in A$ were a branch of
$f_1(t)$ then $a$ would be a cluster point of $B(t)$, and this contradicts with the fact that $B(t)$ is closed and discrete.\\

On the other hand, for every $a\in B(t)$ we define $s_a^t$ to be
the lowest element $s\in\omega^{<\omega}$ such that $s\prec a$ and
$s\not\in f_1(t)$. We define a function
$f_2:\mathcal{T}\To\omega^{(\omega^{<\omega})}$ in the following
way: for $a\in B(t)$, $f_2(t)(s_a^t)$ equals the level of the tree
$T_a$ in which $t$ lies; if $s$ is different from any $s_a^t$,
then say $f_2(t)(s)=0$. Thus, we have defined a function
$f:\mathcal{T}\To wf(A)\times \omega^{(\omega^{<\omega})}$ by
$f(t) = (f_1(t),f_2(t))$, $t\in T$. It
remains to show that this is a determining function.\\

Let $C\subset wf(A)\times \omega^{(\omega^{<\omega})}$ be compact
and we suppose by contradiction that there is an element $x\in
R[\mathcal{F}]\subset 2^\mathcal{T}$, that is a branch $x=\{t_1<_a
t_2<_a \cdots\}$ of the tree $T_a$ for some $a\in A$,
such that $f(x)\subset C$. Two cases arise:\\

Case 1: The elements $s_a^{t_n}$ are equal to some fixed
$s\in\omega^{<\omega}$ for infinitely many $n$'s. For these $n$'s
we have $n\leq f_2(t_n)(s_a^{t_n}) = f_s(t_n)(s)$ and this
contradicts with the fact that $f(x)\subset C$ and $C$ is
compact.\\

Case 2: Modulo passing to a subsequence, we may assume that
$s_a^{t_1}< s_a^{t_2} <\cdots\prec a$. For every $n$ consider the
element $u_n<s_a^{t_n}$ which has length one less than
$s_a^{t_n}$. We have in this case that $u_n\in f_1(t_n)$, and the
$u_n$'s determine that $a$ is a branch of
$\bigcup_{i<\omega}f_1(t_i)$. But, after the Claim below,
$\bigcup_{i<\omega}f_1(t_i)\subset\bigcup_{x\in C}f_1(x)\in
wf(A)$, which
is a contradiction.\\

Claim: If $L\subset wf(A)$ is compact, then $\bigcup L\in wf(A)$.
Proof: Suppose $b\in A$ were a branch of $\bigcup L$. For every
$n<\omega$, let $C_n = \{T\in L : b|_n\in T\}$. This is a
decreasing sequence of nonempty closed subsets of $L$. By
compactness, the intersection of them is nonempty, which implies
that for some $T\in L$, $b$ is a branch of $T$.$\qed$\\

\begin{thm}\label{Rezcompactiscomplex}
Let $A\subset \mathcal{N}$, let $\mathcal{T} = \bigcup_{a\in
A}T_a$, let $\mathcal{F}=\mathcal{F}_A$ be the family of all
closed and discrete subsets of $A$ and let $R[\mathcal{F}]\subset
2^\mathcal{T}$ be the compact set coming from an
$(A,\mathcal{F})$-Reznichenko family of trees. Then, there exists
no determining function $f:\mathcal{T}\To A$.
\end{thm}

For the proof of this result we need Lemma~\ref{bigRez} below,
which is a generalization of~\cite[Lemma 6.2]{ArgArvMer} with
analogous proof. A subset $D$ of a tree $T$ is called
\emph{successively dense} if there is a countable family $R$ of
immediate successors of the root such that every $t\in T$
incomparable with every element of $R$ has an immediate successor
in $D$.

\begin{lem}\label{bigRez}
Let $\{U_n\}_{n<\omega}$ be a disjoint family of open subsets of
$A$ whose union is closed in $A$, and let $\mathcal{T} =
\bigcup_{n<\omega}D_n$ be a countable decomposition of
$\mathcal{T}$. Then, there exists $n<\omega$ such that $D_n\cap
T_a$ is successively dense in the tree $T_a$ for every $a\in U_n$.
\end{lem}

{\tt Proof}: We suppose that the statement of the lemma is false. Then,
we can construct recursively a sequence $(a_n)$ of elements of $A$
and a sequence $(t_n)$ of elements of $\mathcal{T}$ with the
following properties:

\begin{itemize}

 \item $a_n\in U_n$, $t_n\in T_{a_n}$ for every $n$.

\item No immediate successor of $t_n$ in $T_{a_n}$ belongs to
$D_n$.

\item The sets $S_n$ of predecessors of $t_n$ in the tree
$T_{a_n}$ are pairwise disjoint.
\end{itemize}

The construction is performed as follows: Assume that we already defined $a_i$ and $t_i$ for $i<n$. From the negation of the lemma, we obtain $a_n\in U_n$ such that $D_n\cap T_{a_n}$ is not successively dense in $T_{a_n}$. Let $R$ be the set of immediate successors of $a_n$ in $T_{a_n}$ which are comparable with some element of $\bigcup_{i<n}S_i$ in the tree $T_{a_n}$. The set $R$ is finite. Hence, since $D_n\cap T_{a_n}$ is not successively dense in $T_{a_n}$, we can pick $t_n\in T_{a_n}$ incomparable with every $r\in R$ such that no immediate successor of $t_n$ belongs to $D_n$. This finishes the recursive construction.\\

Now, since $a_n\in U_n$ and the $U_n$'s are disjoint open sets
with closed union, the set $\{a_n : n<\omega\}$ is closed and
discrete in $A$, hence it belongs to $\mathcal{F}$. From the
definition of Reznichenko family of trees, we conclude that there
must exist an element $t$ which is simultaneously an immediate
successor of the segment $S_n$ in $T_{a_n}$ for every $n$. For
some $m$, $t\in D_{m}$. But this contradicts that $t_m$ has no
immediate successor in the tree $T_{a_m}$ belonging to
$D_m$.$\qed$\\

{\tt Proof of Theorem~\ref{Rezcompactiscomplex}}: Suppose that there
exists a determining function $f:\mathcal{T}\To A$. Remember that
$A\subset\mathcal{N}$, and that for $s\in\omega^{<\omega}$ we put
$W_s = \{a\in A : s\prec a\}$. Let
$D_s = f^{-1}(W_s)$, and let\\

\begin{itemize}
\item $S_0 = \{s\in\omega^{<\omega} : D_s\text{ is successively
dense in the tree }T_a\text{ for every }a\in W_s\}$.

\item $S_1 = \{s\in S_0 : t\in S_0 \text{ for all } t< s\}$.

\item $S_2 = \{s\in\omega^{<\omega}\setminus S_0 : s\text{ is an
immediate successor of some element of }S_1\}$.
\end{itemize}

Notice that $S_1$ is a tree on $\omega$. We claim that $S_1$ has a
branch $a\in A$. Otherwise no branch of $S_1$ would be an element
of $A$, and this would mean that the union of the family of
disjoint clopen sets $\{W_s : s\in S_2\}$ would be the whole $A$.
We could apply then Lemma~\ref{bigRez} to the decomposition
$\mathcal{T} = \bigcup_{s\in S_2}D_s$ and the clopen sets $\{W_s :
s\in S_2\}$ and we would conclude that there exists $s\in S_2$
such that $D_s$ is successively dense in $T_a$ for every $a\in
W_s$, namely $s\in S_0$ which contradicts that $s\in
S_2$.\\

Let $a\in A$ be a branch of $S_1$. Then $D_{a|n}$ is succesively
dense in $T_a$ for every $n<\omega$. Hence, for every $n<\omega$
there is a countable family $C_n$ of immediate succesors of $a$ in
$T_a$ such that every element of the tree $T_a$ incomparable with
$C_n$ has an immediate succesor in $D_{a|n}$. Let $t$ be an
immediate successor of $a$ in $T_a$ such that
$t\not\in\bigcup_{n<\omega}C_n$. Then, we can construct in the
tree $T_a$ an infinite sequence $t<_a t_1 <_a t_2 <_a\cdots$ with
$t_n\in D_{a|n}= f^{-1}(W_{a|n})$. This infinite sequence
constitues an element of $R[\mathcal{F}]\subset 2^\mathcal{T}$, so
this contradicts that $f$ is a determining function.$\qed$\\

Let us recall the definition of the projective classes $\Sigma^1_n$ and $\Pi^1_n$. As we already indicated $\Sigma^1_1$ and $\Pi^1_1$ denote the classes of analytic and conalytic sets respectively. Recursively, $\Sigma^1_{n+1}$ is defined as the class of separable metrizable spaces which are continuous images of spaces in $\Pi^1_n$, and $\Pi^1_{n+1}$ are the separable metrizable spaces which are complements of sets in $\Sigma^1_{n+1}$ inside a Polish space.\\

{\tt Proof of Theorem~\ref{projectiveclasses}}: Let
$A\subset\omega^\omega\times\omega^\omega$ be a universal
$\Pi^1_{n-1}$ space, that is a $\Pi^1_{n-1}$-set such that for every
$\Pi^1_{n-1}$ subset $B$ of $\omega^\omega$ there exists
$b\in\omega^\omega$ such that $\{b\}\times
B=A\cap(\{b\}\times\omega^\omega)$. Such a set always exists,
cf.~\cite{Kechrisbook}. Set $X=C(R[\mathcal{F}_A])$. By
Theorem~\ref{Rezcompactiscomplex}, for $R[\mathcal{F}_A]\subset
2^\mathcal{T}$, there is no determining function $f:\mathcal{T}\To
A$. \footnote[2]{The rest of this paragraph is different from the published version in Studia Mathematica. It fixes a gap discovered by Marian Fabian.}Since every $\Sigma_n^1$-set is the continuous image of a $\Pi_{n-1}^1$, by Lemma~\ref{determiningcomposition}, it is enough to show that there is no determining function $f:\Gamma\To B$, where $B$ is $\Pi^1_{n-1}$. If such $f$ exists, pick $b$ such that $\{b\}\times
B=A\cap(\{b\}\times\omega^\omega)$. Then the composition $f:\Gamma\To B \sim B\times\{b\} \hookrightarrow A$ would be also a determining function
due to the fact that $B\times\{b\}$ is closed in $A$. This is a
contradiction.\\

Finally, we prove that  $C(R[\mathcal{F}_A])$ is weakly
$\Sigma^1_{n+1}$-determined. By Theorem~\ref{boundcomplexity}, it is enough to show that $wf(A)$ belongs in
$\Pi^1_{n}\subset \Sigma^1_{n+1}$. Put
\[S = \{x\in Tr : \exists a\in A \text{ which is a branch of }x\},\] \[\Omega = \{(x,a) \in Tr\times A : a\text{ is a branch of }x\}.\] Then $\Omega$ is a closed subset of $Tr\times A$, hence a $\Pi^1_{n-1}$-set like $A$. The set $S$ is the projection of $\Omega$ to the first coordinate, so $S$ is $\Sigma^1_n$. Finally, $wf(A) = Tr\setminus S$ is $\Pi^1_n$.$\qed$\\

\section*{Acknowledgements}
I wish to thank Ond\v{r}ej Kalenda for his useful remarks
concerning the work of \v{C}i\v{z}ek and Fabian and also Pandelis
Dodos for teaching to me the descriptive set-theoretic way of
thinking which is somewhat present in this note. An
acknowledgement is also due to the referee, for a number of
corrections and improvements in the final form of the article.

\end{document}